\documentclass[11pt]{article}
\usepackage{epsfig}
\usepackage{amssymb}
\usepackage{latexsym}
\usepackage{epic}
\usepackage{eepic}
\usepackage{amsmath}

\newenvironment{reflist}{\begin{list}{}
         {\itemsep=20pt \parsep=3pt
          \topsep=0pt  \parskip=20pt  \listparindent=-.15in
           \leftmargin= 0.15in }
         \item \ \vspace{-.35in} }
         {\end{list}}

\newcommand{\rb}{\mathbf{r}}
\newcommand{\nb}{\mathbf{n}}
\newcommand{\ep}{\epsilon}
\newcommand{\al}{\alpha}

\newcommand{\beqrn}{ \begin{eqnarray*}}
\newcommand{\eeqrn}{ \end{eqnarray*}}
\newcommand{\beq}{\begin{equation}}
\newcommand{\eeq}{\end{equation}}

\newtheorem{theorem}{Theorem}
\newtheorem{lemma}{Lemma}

\begin{document}

\author{David K\"{a}llberg$^{2}$, Nikolaj Leonenko$^{1}$, Oleg\ Seleznjev$%
^{2}$ \vspace{0.5cm}\hfill \\
$^1 \,$ School of Mathematics, Cardiff University,\\
Senghennydd Road, Cardiff CF24 4YH, UK,\\
$^2 \,$Department of Mathematics and Mathematical Statistics\\
Ume{\aa } University, SE-901 87 Ume\aa , Sweden }
\date{ }
\title{Statistical Inference for R\'enyi Entropy  Functionals}
\maketitle

\begin{abstract}
Numerous entropy-type characteristics (functionals) generalizing R\'enyi
entropy are widely used in mathematical statistics, physics, information
theory, and signal processing for characterizing uncertainty in probability
distributions and distribution identification problems. We consider
estimators of some entropy (integral) functionals for discrete and
continuous distributions based on the number of epsilon-close vector records
in the corresponding independent and identically distributed samples from
two distributions. The estimators form a triangular scheme of generalized $U$%
-statistics. We show the asymptotic properties of these estimators (e.g.,
consistency and asymptotic normality). The results can be applied in various
problems in computer science and mathematical statistics (e.g., approximate
matching for random databases, record linkage, image matching).
\end{abstract}

\baselineskip=3.5 ex



\section{Introduction}

Let $X$ and $Y$ be $d$-dimensional random vectors with discrete or
continuous distributions $\mathcal{P}_{X}$ and $\mathcal{P}_{Y}$,
respectively. In information theory and statistics, there are various
generalizations of Shannon entropy (see Shannon, 1948), characterizing
uncertainty in $\mathcal{P}_{X}$ and $\mathcal{P}_{Y}$, for example, the R%
\'{e}nyi entropy (R\'{e}nyi, 1961, 1970),
\begin{equation*}
h_{s}:=\frac{1}{1-s}\log \left( \int_{R^{d}}p_{X}(x)^{s}dx\right) ,\qquad
s\neq 1,
\end{equation*}%
and the (differentiable) variability for approximate record matching in
random databases
\begin{equation*}
v:=-\log \left( \int_{R^{d}}p_{X}(x)p_{Y}(x)dx\right) ,
\end{equation*}%
where $p_{X}(x),p_{Y}(x),x\in R^{d}$, are densities of $\mathcal{P}_{X}$ and $\mathcal{P}_{Y}$, respectively (see Seleznjev and Thalheim,
2003, 2008). Henceforth we use $\log x$ to denote the natural logarithm of $x$. More generally, for non-negative integers $r_{1},r_{2} \geq
0$ and $\mathbf{r}:=(r_{1},r_{2})$, we consider \emph{R\'{e}nyi entropy functionals}
\begin{equation}
q_{\rb}:=\int_{R^{d}}p_{X}(x)^{r_{1}}p_{Y}(x)^{r_{2}}dx,  \notag
\end{equation}%
and for the discrete case, $\mathcal{P}_{X}=\{p_{X}(k),k\in N^{d}\}$ and $%
\mathcal{P}_{Y}=\{p_{Y}(k),k\in N^{d}\}$,
\begin{equation}
q_{\rb}:=\sum_{k}p_{X}(k)^{r_{1}}p_{Y}(k)^{r_{2}},  \notag
\end{equation}%
i.e., $q_{\rb}=q_{r_{1},r_{2}}$. Then, for example, the R\'{e}nyi entropy $%
h_s=h_{s,0}=\log (q_{s,0})/(1-s)$ and the variability $v=h_{1,1}=-\log (q_{1,1})$. Let $%
X_{1},\ldots ,X_{n_{1}}$ and $Y_{1},\ldots ,Y_{n_{2}}$ be mutually
independent samples of independent and identically distributed (i.i.d.)\
observations from $\mathcal{P}_{X}$ and $\mathcal{P}_{Y}$, respectively. We
consider the problem of estimating the entropy-type functionals $q_{\rb}$
and related characteristics for $\mathcal{P}_{X}$ and $\mathcal{P}_{Y}$ from
the samples $X_{1},\ldots ,X_{n_{1}}$ and $Y_{1},\ldots ,Y_{n_{2}}$.

Various entropy applications in statistics (e.g., classification and
distribution identification problems) and in computer science and
bioinformatics (e.g., average case analysis for random databases,
approximate pattern and image matching) are investigated in, e.g., Kapur
(1989), Kapur and Kesavan (1992), Leonenko et al.\ (2008), Szpankowski
(2001), Seleznjev and Thalheim (2003, 2008), Thalheim (2000), Baryshnikov et
al.\ (2009), and Leonenko and Seleznjev (2010). Some average case analysis
problems for random databases with entropy characteristics are investigated
also in Demetrovics et al.\ (1995, 1998a, 1998b).

In our paper, we generalize the results and approach proposed in Leonenko
and Seleznjev (2010), where the quadratic R\'enyi entropy estimation is
studied for one sample. We consider properties (consistency and asymptotic
normality) of kernel-type estimators based on the number of coincident (or $%
\epsilon$-close) observations in $d$-dimensional samples for more general
class of entropy-type functionals. These results can be used, e.g., in
evaluating of asymptotical confidence intervals for the corresponding R\'{e}nyi
entropy   functionals.

Note that our estimators of entropy-type functionals are different form
those considered by Kozachenko and Leonenko (1987), Tsybakov and van der
Meulen (1996), Leonenko et al.\ (2008), and Baryshnikov et al.\ (2009) (see
Leonenko and Seleznjev, 2010, for a discussion).

First we introduce some notation. Throughout the paper, let $X$ and $Y$ be
independent random vectors in $R^{d}$ with distributions $\mathcal{P}_{X}$
and $\mathcal{P}_{Y}$, respectively. For the discrete case, $\mathcal{P}%
_{X}=\{p_{X}(k),k\in N^{d}\}$ and $\mathcal{P}_{Y}=\{p_{Y}(k),k\in N^{d}\}$.
In the continuous case, let the distributions be with densities $p_{X}(x)$
and $p_{Y}(x),x\in R^{d}$, respectively. Let $d(x,y)=|x-y|$ denote the
Euclidean distance in $R^{d}$ and $B_{\ep}(x):=\{y:d(x,y)\leq \epsilon \}$
an $\epsilon $-ball in $R^{d}$ with center at $x$, radius $\epsilon $, and
volume $b_{\ep}(d)=\epsilon ^{d}b_{1}(d)$, $b_{1}(d)=2\pi ^{{d}/{2}%
}/(d\Gamma( {d}/{2}) )$. Denote by $p_{X,\epsilon }(x)$ the $%
\epsilon $-ball probability
\begin{equation}
p_{X,\epsilon }(x):=P\{X\in B_{\ep}(x)\}.  \notag
\end{equation}%
We write $I(C)$ for the indicator of an event $C$, and $|D|$ for the
cardinality of a finite set $D$.

Next we define the following estimators of $q_{\mathbf{r}}$ when $r_1$ and $%
r_2$ are non-negative integers. Let the i.i.d.\ samples $X_1,\ldots,X_{n_1}$
and $Y_1,\ldots,Y_{n_2}$ be from $\mathcal{P}_X$ and $\mathcal{P}_Y$,
respectively. Denote $\mathbf{n}:= (n_1,n_2)$, $n:=n_1+n_2$, and say that $%
\mathbf{n} \to \infty$ if $n_1, n_2 \to \infty$ and let $p_\nb := n_1/n \to
p, 0<p<1$, as $\mathbf{n} \to \infty$.

For an integer $k$, denote by $\mathcal{S}_{m,k}$ the set of all $k$-subsets
of $\{1,\ldots,m\}$. For $S \in \mathcal{S}_{n_1,r_1}$, $T \in \mathcal{S}%
_{n_2,r_2}$, and $i \in S$, define
\begin{equation}
\psi_\nb^{(i)}(S;T) := I(d(X_i,X_j)\leq \epsilon,d(X_i,Y_k)\leq
\epsilon,\forall j \in S, \forall k \in T),  \notag
\end{equation}
i.e., the indicator of the event that all elements in $\{X_j, j \in S\}$ and
$\{Y_k, k \in T\}$ are $\epsilon$-close to $X_i$. Note that by conditioning
we have
\begin{equation*}
\mathrm{E} \psi_\nb^{(i)}(S;T)= \mathrm{E} p_{X,\epsilon}(X)^{r_1-1}p_{Y,%
\epsilon}(X)^{r_2} =: q_{\mathbf{r},\epsilon} ,
\end{equation*}
say, the $\epsilon$-coincidence probability. Let a generalized $U$-statistic
for the functional $q_{\mathbf{r},\epsilon}$ be defined as
\begin{eqnarray*}  \label{statistic}
Q_\nb=Q_{\mathbf{n},\mathbf{r},\epsilon}:= \binom{n_1}{r_1}^{-1}\binom{n_2}{%
r_2}^{-1}\sum_{(n_1,r_1)}\sum_{(n_2,r_2)} \psi_\nb(S;T),  \notag \\
\end{eqnarray*}
where the symmetrized kernel
\begin{equation*}
\psi_\nb(S;T) := \frac{1}{r_1}\sum_{i \in S}\psi_\nb^{(i)}(S;T),
\end{equation*}
and by definition, $Q_\nb$ is an unbiased estimator of $q_{\mathbf{r}%
,\epsilon}= \mathrm{E} Q_\nb$. Define for discrete and continuous
distributions
\begin{eqnarray*}
\zeta_{1,0} &:=& \mathrm{Var}%
(p_X(X)^{r_1-1}p_Y(X)^{r_2})=q_{2r_1-1,2r_2}-q_{r_1,r_2}^2,  \notag \\
\zeta_{0,1} &:=& \mathrm{Var}%
(p_X(Y)^{r_1}p_Y(Y)^{r_2-1})=q_{2r_1,2r_2-1}-q_{r_1,r_2}^2,  \notag \\
\kappa &:=& p^{-1}r_1^2 \zeta_{1,0}+(1-p)^{-1}r_2^2 \zeta_{0,1}.  \notag
\end{eqnarray*}
Let $\overset{\mathrm{D}}{\rightarrow }$ and $\overset{\mathrm{P}}{%
\rightarrow }$ denote convergence in distribution and in probability,
respectively.

The paper is organized as follows. In Section \ref{se:Main}, we consider
estimation of R\'enyi entropy  functionals for discrete and continuous
distributions. In Section \ref{se:App}, we discuss some applications of the
obtained estimators in average case analysis for random databases (e.g., for
join optimization with approximate matching), in pattern and image matching
problems, and for some distribution identification problems. Several numerical
experiments demonstrate the rate of convergence in the obtained asymptotic
results. Section \ref{se:Proofs} contains the proofs of the statements from
the previous sections.


\section{Main Results}

\label{se:Main}

\subsection{Discrete Distributions}

In the discrete case, set $\epsilon = 0$, i.e., exact coincidences are
considered. Then $Q_\nb$ is an unbiased estimator of the $\epsilon$%
-coincidence probability
\begin{equation*}
q_{\mathbf{r},0}=q_\rb = \mathrm{E} I(X_1=X_i=Y_j, i = 2,\ldots, r_1, j =
1,\ldots, r_2) = \mathrm{E} p_X(X)^{r_1-1}p_Y(X)^{r_2}.
\end{equation*}
Let $Q_{\mathbf{n},\mathbf{r}} :=Q_{\mathbf{n},\mathbf{r},0}$ and
\begin{equation*}
K_\nb:=p_\nb^{-1} r_1^2 (Q_{\mathbf{n},2r_1-1,2r_2}-Q_{\mathbf{n},\mathbf{r}%
}^2) + (1-p_\nb)^{-1}r_2^2(Q_{\mathbf{n},2r_1,2r_2-1}-Q_{\mathbf{n},\mathbf{r%
}}^2),
\end{equation*}
and $k_\nb:=\max(K_\nb,1/n)$, an estimator of $\kappa$. Denote by $H_\nb :=
\log(\max(Q_\nb,1/n))/(1-r)$, an estimator of $h_\rb := \log(q_\rb)/(1-r)$.
\newline
\newline
\textbf{Remark.} Instead of $1/n$ in the definition of a truncated estimator, a sequence $a_n>0,$ $a_n \to 0$ as $\mbox{ as } \mathbf{n}
\to \infty$, can be used (cf.\ Leonenko and Seleznjev, 2010).

The next asymptotic normality theorem for the estimator $Q_\nb$ follows
straightforwardly from the general $U$-statistic theory (see, e.g., Lee,
1990, Koroljuk and Borovskich, 1994) and the Slutsky theorem.

\begin{theorem}
\label{th:Disc} If $\zeta_{1,0},\zeta_{0,1} > 0$, then
\begin{align}
& \sqrt{n}(Q_\nb-q_\rb) \overset{\mathrm{D}}{\to} N(0, \kappa) \mbox{ and }
\sqrt{n}(Q_\nb-q_\rb)/k_\nb^{1/2} \overset{\mathrm{D}}{\to} N(0,1);  \notag
\\
& \sqrt{n}(1-r)\frac{Q_\nb}{k_\nb^{1/2}}(H_\nb-h_\rb) \overset{\mathrm{D}}{%
\to} N(0,1) \mbox{ as } \mathbf{n} \to \infty.  \notag
\end{align}
\end{theorem}

\subsection{Continuous Distributions}

In the continuous case, denote by $\tilde{Q}_\nb := Q_{\mathbf{n}}/b_\ep
(d)^{r-1}$ an estimator of $q_{\mathbf{r}}$. Let $\tilde{q}_{\mathbf{r}%
,\epsilon} := \mathrm{E}\tilde{Q}_\nb=q_{\mathbf{r},\epsilon}/b_\ep
(d)^{r-1} $ and $v^2_\nb:= \mathrm{Var}(\tilde{Q}_\nb)$.

Henceforth, assume that $\epsilon = \epsilon(\mathbf{n}) \to 0$ as $\mathbf{n%
} \to \infty$. For a sequence of random variables $U_n, n \geq 1$, we say
that $U_n = \mathrm{O}_{\mathrm{P}}(1)$ as $n \to \infty$ if for any $%
\epsilon > 0$ and $n$ large enough there exists $A>0$ such that $P(|U_n|>A)
\leq \epsilon$, i.e., the family of distributions of $U_n, n\geq 1$, is
tight, and for a numerical sequence $w_n, n \geq 1$, say, $U_n = \mathrm{O}_{%
\mathrm{P}}(w_n)$ as $n \to \infty$ if $U_n/w_n = \mathrm{O}_{\mathrm{P}}(1)$
as $n \to \infty$. The following theorem describes the consistency and
asymptotic normality properties of the estimator $\tilde{Q}_\nb$.

\begin{theorem}
\label{thm:main} Let $p_X(x)$ and $p_Y(x)$ be bounded and continuous or with
a finite number of discontinuity points.
\itemize

\item[(i)] $v^2_\nb = \mathrm{O}(n^{-1}\epsilon^{d(1/r-1)})$ and $\mathrm{E}
\tilde{Q}_\nb \to q_\rb$ as $\mathbf{n} \to \infty$, and hence if $%
n\epsilon^{d(1-1/r)} \to \infty$ as $\mathbf{n} \to \infty$, then $\tilde{Q}%
_\nb$ is a consistent estimator of $q_\rb$.

\item[(ii)] If $n\epsilon^d \to \infty$ as $\mathbf{n} \to \infty$ and $%
\zeta_{1,0},\zeta_{0,1} > 0$, then
\begin{eqnarray*}  \label{main}
\sqrt{n}(\tilde{Q}_\nb-\tilde{q}_{\mathbf{r},\epsilon}) &\overset{\mathrm{D}}%
{\to} & N(0,\kappa)\; \mbox{ as } \mathbf{n} \to \infty.  \notag
\end{eqnarray*}
\end{theorem}

\bigskip

In order to evaluate the functional $q_\rb$, we denote by $H^{(\alpha)}(C)$ $%
, 0<\alpha\le 2,C>0 $, a linear space of bounded and continuous in $R^d$
functions satisfying $\alpha$-H\"{o}lder condition if $0<\alpha\le 1$ or if $%
1<\alpha\le 2$ with continuous partial derivatives satisfying $(\alpha-1)$-H%
\"{o}lder condition with constant $C$. Furthermore, let
\begin{equation}
K_\nb:= p_\nb^{-1} r_1^2 (\tilde{Q}_{\mathbf{n},2r_1-1,2r_2,\epsilon}-\tilde{%
Q}_{\mathbf{n},\mathbf{r},\epsilon}^2) + (1-p_\nb)^{-1}r_2^2(\tilde{Q}_{%
\mathbf{n},2r_1,2r_2-1,\epsilon}-\tilde{Q}_{\mathbf{n},\mathbf{r}%
,\epsilon}^2),  \notag
\end{equation}
and define $k_\nb := \max(K_\nb,1/n)$. It follows from Theorem \ref{thm:main}
and Slutsky's theorem that $k_\nb$ is a consistent estimator of the
asymptotic variance $\kappa$. Denote by $H_\nb := \log(\max(\tilde{Q}%
_\nb,1/n))/(1-r)$, an estimator of $h_\rb := \log(q_\rb)/(1-r)$. Let $L(n)$
be a slowly varying function. We obtain the following asymptotic result.

\begin{theorem}
\label{th:AsBias} Let $p_X(x), p_Y(x) \in H^{(\alpha)}(C)$.\newline
(i) Then the bias $\; |\tilde{q}_{\mathbf{r},\epsilon}-q_\rb|\le
C_1\epsilon^{\alpha}, C_1>0$.\newline
(ii) If $0<\alpha\le d/2$ and $\epsilon\sim c
n^{-\alpha/(2\alpha+d(1-1/r))}, 0<c<\infty$, then
\begin{eqnarray*}
\tilde{Q}_\nb-q_\rb = \mathrm{O}_{\mathrm{P}}(n^{-\alpha/(2%
\alpha+d(1-1/r))}) \mbox{  and  } {H_\nb}-h_\rb = \mathrm{O}_{\mathrm{P}%
}(n^{-\alpha/(2\alpha+d(1-1/r))}) \mbox{ as } \mathbf{n} \to \infty.
\end{eqnarray*}
(iii) If $\alpha > d/2$ and $\epsilon \sim L(n)n^{-1/d}$ and $n\epsilon^d
\to \infty$, then
\begin{align}
& \sqrt{n}(\tilde{Q}_\nb-q_\rb) \overset{\mathrm{D}}{\to} N(0,\kappa)
\mbox{
and } \sqrt{n}(\tilde{Q}_\nb-q_\rb)/k_\nb^{1/2} \overset{\mathrm{D}}{\to}
N(0,1);  \notag \\
& \sqrt{n} (1-r)\frac{\tilde{Q}_\nb}{k_\nb^{1/2}}(H_\nb-h_\rb) \overset{%
\mathrm{D}}{\to} N(0,1) \mbox{ as } \mathbf{n} \to \infty.  \notag
\end{align}
\end{theorem}

\section{Applications and Numerical Experiments}

\label{se:App}

\subsection{Approximate Matching in Stochastic Databases}

Let tables (in a \emph{relational database}) $T_1$ and $T_2$ be matrices
with $m_1$ and $m_2$ i.i.d.\ random tuples (or records), respectively. One
of basic database operations, \emph{join}, combines two tables into a third
one by matching values for given columns (attributes). For example, the join
condition can be the equality (equi-join) between a given pairs of
attributes (e.g., names) from the tables. Joins are especially important for
tieing together pieces of disparate information scattered throughout a
database (see, e.g., Kiefer et al.\ 2005, Copas and Hilton, 1990, and
references therein). For the approximate join, we match $\epsilon$-close
tuples, say, $d(t_1(j),t_2(i))\le \epsilon, t_k(j)\in T_k, k=1,2$, with a
specified distance, see, e.g., Seleznjev and Thalheim (2008). A set of
attributes $A$ in a table $T$ is called an $\epsilon$-key (test) if there
are no $\epsilon$-close sub-tuples $t_A(j), j=1,\dots,m$. Knowledge about
the set of tests ($\epsilon$-keys) is very helpful for avoiding redundancy
in identification and searching problems, characterizing the complexity of a
database design for further optimization, see, e.g., Thalheim (2000). By
joining a table with itself (self-join) we identify also $\epsilon$-keys and
key-properties for a set of attributes for a random table (Seleznjev and
Thalheim, 2003, Leonenko and Seleznjev, 2010).

The cost of join operations is usually proportional to the size of the
intermediate results and so the joining order is a primary target of
join-optimizers for multiple (large) tables, Thalheim (2000). The average
case approach based on stochastic database modelling for optimization
problems is proposed in Seleznjev and Thalheim (2008), where for random
databases, the distribution of the $\epsilon$-join size $N_\ep$ is studied.
In particular, with some conditions it is shown that the average size
\begin{equation*}
\mathrm{E} N_\ep = m_1 m_2 q_{1,1,\epsilon}= m_1 m_2 \epsilon^d b_1(d)
(e^{-h_{1,1}} + \mathrm{o}(1)) \mbox{ as } \epsilon\to 0,
\end{equation*}
that is the asymptotically optimal (in average) pairs of tables are amongst
the tables with maximal value of the functional $h_{1,1}$ (variability) and
the corresponding estimators of $h_{1,1}$ can be used for samples $%
X_1,\ldots, X_{n_1}$ and $Y_1,\ldots, Y_{n_2}$ from $T_1$ and $T_2$,
respectively. For discrete distributions, similar results from Theorem \ref%
{th:Disc} for $\epsilon=0$ can be applied.


\subsection{\textbf{Image Matching using Quadratic-entropy Measures}}

%
%

Image retrieval and registration fall in the general area of pattern
matching problems, where the best match to a reference or query image $I_{0}$
is to be found in a database of secondary images $\left\{ I_{i}\right\}
_{i=1}^{n}.$ The best match is expressed as a partial re-indexing of the
database in decreasing order of similarity to the reference image using a
similarity measure. In the context of image registration, the database
corresponds to an infinite set of transformed versions of a secondary image,
e.g., rotation and translation, which are compared to the reference image to
register the secondary one to the reference.

Let $X$ be a $d$-dimensional random vector and let $p(x)$ and $q(x)$ denote
two possible densities for $X$. In the sequel, $X$ is a feature vector
constructed from the query image and a secondary image in an image database
and $p(x)$ and $q(x)$ are densities, e.g., for the query image features and
the secondary image features, respectively, say, image densities. When the
features are discrete valued the $p(x)$ and $q(x)$ are probability mass
functions.

The basis for entropy methods of image matching is a measure of similarity
between image densities. A general entropy similarity measure is the R\'{e}%
nyi $\alpha $-divergence, also called the R\'{e}nyi $\alpha $-relative
entropy, between $p(x)$ and $q(x)$%
\begin{equation*}
D_{\alpha}(p,q)=\frac{1}{\alpha-1}\log \int_{R^{d}}q(x)\left( \frac{p(x)}{%
q(x)}\right) ^{\alpha}dx=\frac{1}{\alpha-1}\log
\int_{R^{d}}p^{\alpha}(x)q^{1-\alpha}(x)dx,\;\alpha \neq 1.
\end{equation*}%
When the density $p(x)$ is supported on a compact domain and $q(x)$ is
uniform over this domain, the R\'{e}nyi $\alpha$-divergence reduces to the R%
\'{e}nyi $\alpha$-entropy
\begin{equation*}
h_{\alpha}(p)=\frac{1}{1-\alpha}\log \int_{R^d} p^{\alpha}(x)dx.
\end{equation*}

%
Another important example of statistical distance between distributions is
given by the following nonsymmetric Bregman distance (see, e.g.,  Pardo, 2006)
\begin{equation*}
B_s(p,q)= \int_{R^d} \left [ q(x)^s + \frac{1}{s-1}p(x)^s-\frac{s}{s-1}%
p(x)q(x)^{s-1} \right]dx, \qquad s \neq 1,
\end{equation*}
or its symmetrized version
\begin{equation*}
K_s(p,q) = \frac{1}{s}[B_s(p,q)+B_s(q,p)]= \frac{1}{s-1}\int_{R^d}
[p(x)-q(x)][p(x)^{s-1}-q(x)^{s-1}]dx.
\end{equation*}
For $s = 2$, we get the second order distance
\begin{equation*}
B_2(p,q)=K_2(p,q)= \int_{R^d} [p(x)-q(x)]^2 dx.
\end{equation*}
Now, for an integer $s$, applying Theorem \ref{main} and \ref{th:AsBias} one
can obtain an asymptotically normal estimator of the R\'{e}nyi $s$-entropy
and a consistent estimator of the Bregman distance.

\subsection{Entropy Maximizing Distributions}

For a positive definite and symmetric matrix $\Sigma$, $s \neq 1$, define
the constants
\begin{equation*}
m = d + 2/(s-1), \qquad \mathbf{C}_s = (m+2)\Sigma,
\end{equation*}
and
\begin{equation*}
A_s = \frac{1}{|\pi \mathbf{C}_s|^{1/2}}\frac{\Gamma(m/2+1)}{%
\Gamma((m-d)/2+1)}.
\end{equation*}
Among all densities with mean $\mu$ and covariance matrix $\Sigma$, the R%
\'{e}nyi entropy $h_s$, $s = 2, \ldots,$ is uniquely maximized by the
density (Costa et al.\ 2003)
\begin{equation}  \label{density}
p_s^*(x)= \left \{
\begin{array}{lr}
A_s(1-(x-\mu)^T \mathbf{C}_s^{-1}(x-\mu))^{1/(s-1)}, & x \in \Omega_s \\
0, & x \notin \Omega_s,%
\end{array}
\right.
\end{equation}
with support
\begin{equation*}
\Omega_s = \{x \in R^d:(x-\mu)^T\mathbf{C}^{-1}_s (x-\mu) \leq 1\}.
\end{equation*}
The distribution given by $p_s^*(x)$ belongs to the class of \textit{%
Student-r} distributions. Let $\mathcal{K}$ be a class of $d$-dimensional
density functions $p(x)$, $x \in R^d$, with positive definite covariance
matrix. By the procedure described in Leonenko and Seleznjev (2010), the
proposed estimator of $h_s$ can be used for distribution identification
problems, i.e., to test the null hypothesis $H_0: X_1, \ldots, X_n$ \textit{%
is a sample from a Student-r distribution of type \eqref{density}} against
the alternative $H_1: X_1, \ldots, X_n$ \textit{is a sample from any other
member of $\mathcal{K}$}.

\subsection{Numerical Experiments}

\label{se:Num}

\begin{figure}[hbt]
\begin{center}
\includegraphics[width=0.85\textwidth,height=0.3\textheight]{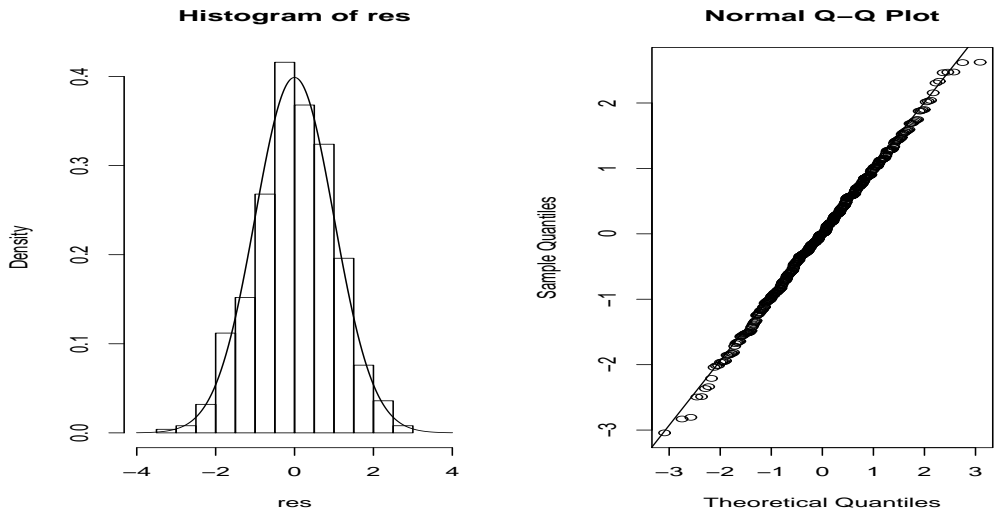}
\end{center}
\caption{Bernoulli $d$-dimensional distribution; $d=3$, $Be(p)$-i.i.d.
components, $p=0.8$, sample size $n=200$. Standard normal approximation for
the empirical distribution (histogram) for the normalized residuals, $%
N_{sim}= 500$.}
\label{fg:Ex1}
\end{figure}

\textbf{Example 1.} Figure \ref{fg:Ex1} shows the accuracy of the estimator
for the cubic R\'{e}nyi entropy $h_{3,0}$ of discrete distributions in
Theorem \ref{th:Disc}, for a sample from a $d$-dimensional Bernoulli
distribution and $n$ observations, $d=3$, $n=300,$ with Bernoulli $Be(p)$%
-i.i.d. components, $p=0.8$. Here the coincidence probability $q_{3,0} =
(p^3+(1-p)^3)^3$ and the R\'{e}nyi entropy $h_{3,0} = -\log(q_{3,0})/2$. The
histogram for the normalized residuals $r^{(i)}_n := 2\sqrt{n}%
Q_\nb(H_\nb-h_\rb)/k_\nb^{1/2},$ $i=1,\ldots, N_{sim}$ are compared to the
standard normal density, $N_{sim} = 500$. The corresponding qq-plot and
p-values for the Kolmogorov-Smirnov (0.4948) and Shapiro-Wilk (0.7292) tests
also support normality hypothesis for the obtained residuals.

\bigskip

\begin{figure}[hbt]
\begin{center}
\includegraphics[width=0.85\textwidth,height=0.3\textheight]{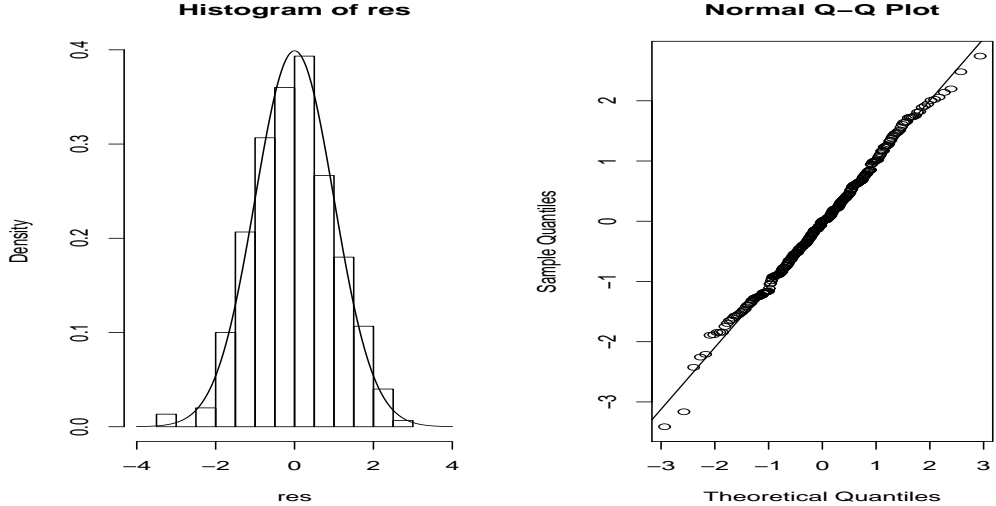}
\end{center}
\caption{Two Gaussian distributions; $N(0,3/2)$, $N(2,1/2)$, $n_1=100, n_2 =
200, \protect\epsilon = 1/10$. Standard normal approximation for the
empirical distribution (histogram) for the normalized residuals, $N_{sim}=
300$.}
\label{fg:Ex2}
\end{figure}

\noindent \textbf{Example 2.} Figure \ref{fg:Ex2} illustrates the
performance of the approximation for the differentiable variability $%
v=h_{1,1}$ in Theorem \ref{th:AsBias}, for two one-dimensional samples from
normal distributions $N(0,3/2)$ and $N(2,1/2)$, with the sample sizes $n_1 =
100,n_2 = 200$, respectively. Here the variability $v=\log(2\sqrt{\pi}e)$.
The normalized residuals are compared to the standard normal density, $%
N_{sim}=300$. The qq-plot and p-values for the Kolmogorov-Smirnov (0.9916)
and Shapiro-Wilk (0.5183) tests also support the normal approximation.

\bigskip

\begin{figure}[hbt]
\begin{center}
\includegraphics[width=0.85\textwidth,height=0.3\textheight]{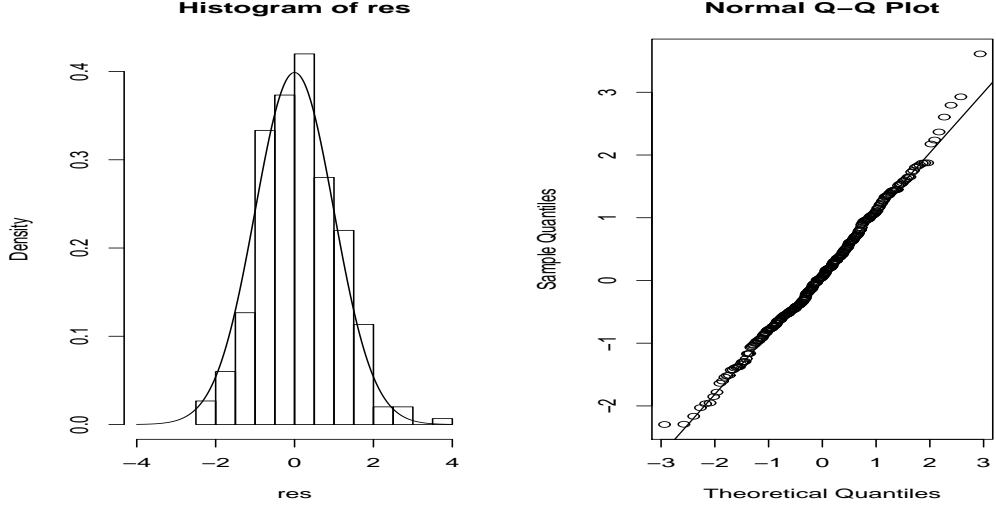}
\end{center}
\caption{Bivariate normal distribution with $N(0,1)$-i.i.d. components;
sample size $n=300, \protect\epsilon = 1/2$. Standard normal approximation
for the empirical distribution (histogram) for the normalized residuals, $%
N_{sim}= 300$.}
\label{fg:Ex3}
\end{figure}

\noindent \textbf{Example 3.} Figure \ref{fg:Ex3} shows the accuracy of the
normal approximation for the cubic R\'{e}nyi entropy $h_{3,0}$ in Theorem %
\ref{th:AsBias}, for a sample from a bivariate Gaussian distribution with $%
N(0,1)$-i.i.d. components, and $n=300$ observations. Here the R\'{e}nyi
entropy $h_{3,0} = \log(\sqrt{12}\pi)$. The histogram, qq-plot, and p-values
for the Kolmogorov-Smirnov (0.2107) and Shapiro-Wilk (0.2868) tests allow to
accept the hypothesis of standard normality for the residuals, $N_{sim}= 300$%
.

\bigskip

\begin{figure}[hbt]
\begin{center}
\includegraphics[width=8.5cm]{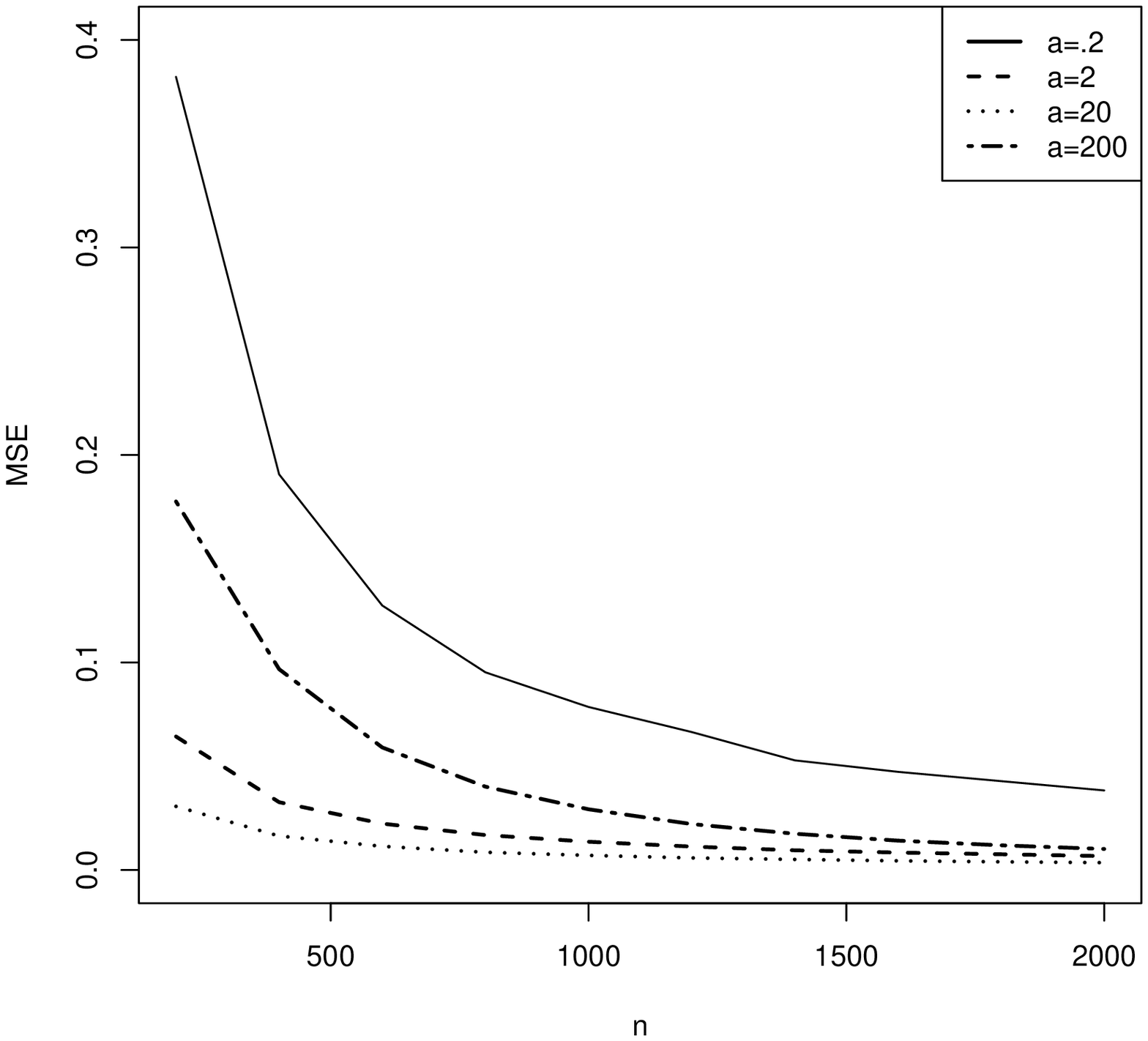}
\end{center}
\caption{Bregman distance for $Exp(\protect\beta_1)$ and $Exp(\protect\beta%
_2)$, $\protect\beta_1 = 1, \protect\beta_2 = 3$. The empirical MSE obtained
for the $U$-statistic estimator with $n\protect\epsilon = a$, for different
values of a.}
\label{fg:Ex4}
\end{figure}

\noindent \textbf{Example 4.} Figure \ref{fg:Ex4} demonstrates the behaviour
of the estimator for the quadratic Bregman distance $B_2(p,q)$ for two
exponential distributions $p(x) = \beta_1e^{-\beta_1 x},x > 0$, and $q(x) =
\beta_2e^{-\beta_2 x}, x >0$, with rate parameters $\beta_1 = 1, \beta_2 = 3$%
, respectively, and equal sample sizes. Here the Bregman distance $B_2(p,q)
= 1/2$. The empirical mean squared error (MSE) based on 10000 independent
simulations are calculated for different values of $n$.


%
%
%
%
%
%
%

\section{Proofs}

\label{se:Proofs}

\begin{lemma}
Assume that $p_X(x)$ and $p_Y(x)$ are bounded and continuous or with a
finite number of discontinuity points. Let $a,b \geq 0$. Then
\begin{equation*}
b_\ep (d)^{-(a+b)}\mathrm{E}(p_{X,\epsilon}(X)^a p_{Y,\epsilon}(X)^b) \to
\int_{R^d}p_X(x)^{a+1} p_Y(x)^{b} dx \mbox{ as } \epsilon \to 0.  \notag
\end{equation*}
\end{lemma}

\textit{Proof: } We have
\begin{equation}
b_\ep (d)^{-(a+b)}\mathrm{E}(p_{X,\epsilon}(X)^{a}p_{Y,\epsilon}(X)^{b}) =
\mathrm{E}(g_\ep(X)),  \notag
\end{equation}
where $g_\ep(x):=(p_{X,\epsilon}(x)/b_\ep (d))^{a}(p_{Y,\epsilon}(x)/b_\ep
(d))^{b}$. It follows by definition that $g_\ep(x) \to p_X(x)^ap_Y(x)^b$ as $%
\epsilon \to 0$, for all continuity points of $p_X(x)$ and $p_Y(x)$, and
that the random variable $g_\ep(X)$ is bounded. Hence, the bounded
convergence theorem implies
\begin{equation}
\mathrm{E}(g_\ep(X)) \to \mathrm{E}(p_X(X)^{a}p_Y(X)^b) = q_{a+1,b}
\mbox{
as } \epsilon \to 0.  \notag
\end{equation}
\hfill $\Box$\newline
\noindent \textit{Proof of Theorem \ref{thm:main}: }$(i)$ Note that for $k
=1,\ldots,r$,
\begin{equation}  \label{eq:k}
n^k \epsilon^{d(k-1)} = (n \epsilon^{d(1-1/k)})^k \geq (n
\epsilon^{d(1-1/r)})^{k} \geq n \epsilon^{d(1-1/r)}.
\end{equation}
We use the conventional results from the theory of $U$-statistics (see,
e.g., Lee, 1990, Koroljuk and Borovskich, 1994). For $l = 0,\ldots, r_1,$
and $m = 0, \ldots, r_2$, define
\begin{align}  \label{cond}
\psi_{l,m,\mathbf{n}}(x_1, &\ldots,x_l;y_1,\ldots,y_m) := \mathrm{E}%
\psi_\nb(x_1,\ldots,x_l,X_{l+1},\ldots,X_{r_1};y_1,\ldots,y_m,Y_{m+1},%
\ldots,Y_{r_2})  \notag \\
&= \frac{1}{r_1}\sum_{i=1}^{r_1}\mathrm{E} \psi_\nb^{(i)}(x_1,%
\ldots,x_l,X_{l+1},\ldots,X_{r_1};y_1,\ldots,y_m,Y_{m+1},\ldots,Y_{r_2}),
\end{align}
and
\begin{equation*}
\sigma^2_{{l},{m},\epsilon} :=\mathrm{Var}(\psi_{l,m,\mathbf{n}%
}(X_1,\ldots,X_l;Y_1,\ldots,Y_m)).
\end{equation*}
Let $S_1,S_2 \in \mathcal{S}_{n_1,r_1}$ and $T_1,T_2 \in \mathcal{S}%
_{n_2,r_2}$ have ${l}$ and ${m}$ elements in common, respectively. By
properties of $U$-statistics, we have 
\begin{equation}  \label{eq:var}
v^2_\nb = \mathrm{Var}(\tilde{Q}_\nb) = b_\ep (d)^{-2(r-1)}\sum_{{l}%
=0}^{r_1}\sum_{{m}=0}^{r_2}\frac{\binom{r_1}{{l}}\binom{r_2}{{m}}\binom{%
n_1-r_1}{r_1-l}\binom{n_2-r_2}{r_2-{m}}} {\binom{n_1}{r_1}\binom{n_2}{r_2}}
\sigma^2_{{l},{m},\epsilon},
\end{equation}
and
\begin{equation}  \label{eq:Sigma}
\sigma^2_{{l},{m},\epsilon}= \mathrm{Cov}(\psi_\nb(S_1;T_1),\psi_%
\nb(S_2;T_2)).
\end{equation}
From \eqref{eq:Sigma} we get that $0 \leq \sigma^2_{{l},{m},\epsilon} \leq
\mathrm{E}(\psi_\nb(S_1;T_1)\psi_\nb(S_2;T_2))$, which is a finite linear
combination of {$P(A_i \cap A_j), i \in S_1, j \in S_2$}, where
\begin{equation}  \label{eq:tri1}
A_i := \{d(X_i,X_k) \leq \epsilon,d(X_i,Y_s)\leq \epsilon, \forall k \in
S_1,\forall s \in T_1\}.  \notag
\end{equation}
When $l \neq 0$ or $m \neq 0$, the triangle inequality implies that
\begin{equation}
A_i \cap A_j \subseteq F_i := \{ d(X_i,X_k)\leq 3\epsilon,d(X_i,Y_s) \leq
3\epsilon,\forall k \in S_1 \cup S_2 ,\forall s \in T_1 \cup T_2\},  \notag
\end{equation}
and since $|S_1 \cup S_2| = 2r_1-l$ and $|T_1 \cup T_2| = 2r_2-m$, it
follows by conditioning and from Lemma 1 that
\begin{eqnarray*}
P(A_i \cap A_j) \leq P(F_i) &= &\mathrm{E}(p_{X,3%
\epsilon}(X_i)^{2r_1-l-1}p_{Y,3\epsilon}(X_i)^{2r_2-m})  \notag \\
&\sim & 3^{d(2r-{l}-{m}-1)} b_\ep (d)^{2r-{l}-{m}-1}q_{2r_1-l,2r_2-m}
\mbox{
as } \mathbf{n} \to \infty.  \notag
\end{eqnarray*}
We conclude that
\begin{equation}  \label{eq:bound}
\sigma^2_{{l},{m},\epsilon} = \mathrm{O}(b_\ep (d)^{2r-{l}-{m}-1})
\mbox{ as
} \mathbf{n} \to \infty.
\end{equation}
Now, for $l=0,\ldots,r_1$ and $m=0,\ldots,r_2,$ we obtain
\begin{equation}  \label{eq:varlim}
b_\ep (d)^{-2(r-1)} \frac{\binom{r_1}{{l}}\binom{r_2}{{m}}\binom{n_1-r_1}{%
r_1-l}\binom{n_2-r_2}{r_2-{m}}} {\binom{n_1}{r_1}\binom{n_2}{r_2}} \sigma^2_{%
{l},{m},\epsilon} \sim C_{l,m} \frac{b_\ep (d)^{-(2r-l-m-1)}\sigma^2_{{l},{m}%
,\epsilon}}{n^{l+m}\epsilon^{d(l+m-1)}} \mbox{ as } \mathbf{n} \to \infty,
\end{equation}
for some constant $C_{l,m}> 0$. Hence, from \eqref{eq:k}, \eqref{eq:var}, %
\eqref{eq:bound}, and \eqref{eq:varlim} we get that $v^2_\nb = \mathrm{O}
((n\epsilon^{d(1-1/r)})^{-1})$ as $\mathbf{n} \to \infty$. Moreover, it
follows from Lemma 1 that $\mathrm{E} \tilde{Q}_\nb \to q_\rb$ as $\mathbf{n}
\to \infty$, so when $n\epsilon^{d(1-1/r)} \to \infty$, then
\begin{equation*}
\mathrm{E}(\tilde{Q}_\nb-q_\rb)^2=v^2_\nb + (\mathrm{E}(\tilde{Q}_\nb-q_{%
\mathbf{r}}))^2 \to 0,
\end{equation*}
and the assertion follows. \newline
\newline
\textit{(ii)} Let
\begin{equation}
h_\nb^{(1,0)}(x) := \psi_{1,0,\mathbf{n}}(x)/b_\ep (d)^{r-1}-\tilde{q}_{%
\mathbf{r},\epsilon}, \quad h_\nb^{(0,1)}(x) := \psi_{0,1,\mathbf{n}%
}(x)/b_\ep (d)^{r-1}-\tilde{q}_{\mathbf{r},\epsilon}.
\end{equation}
The H-decomposition of $\tilde{Q}_\nb$ is given by
\begin{equation}  \label{eq:dec}
\tilde{Q}_\nb = \tilde{q}_{\mathbf{r},\epsilon}+r_1 H_\nb^{(1,0)}+r_2
H_\nb^{(0,1)}+ R_\nb,
\end{equation}
where
\begin{eqnarray*}
H_\nb^{(1,0)} := \frac{1}{n_1}\sum_{i=1}^{n_1}h_\nb^{(1,0)}(X_i), \quad
H_\nb^{(0,1)} := \frac{1}{n_2}\sum_{i=1}^{n_2}h_\nb^{(0,1)}(Y_i).  \notag
\end{eqnarray*}
The terms in \eqref{eq:dec} are uncorrelated, and since $\mathrm{Var}%
(h_\nb^{(1,0)}(X_1)) = b_\ep (d)^{-2(r-1)}\sigma^2_{1,0,\epsilon}$ and $%
\mathrm{Var}(h_\nb^{(0,1)}(Y_1)) = b_\ep
(d)^{-2(r-1)}\sigma^2_{0,1,\epsilon} $, we obtain from \eqref{eq:var} that
\begin{eqnarray}  \label{eq:varrem}
\mathrm{Var}(R_\nb) &=& \mathrm{Var}(\tilde{Q}_\nb) - \mathrm{Var}(r_1
H_\nb^{(1,0)})-\mathrm{Var}(r_2 H_\nb^{(0,1)})  \notag \\
&=&\mathrm{Var}(\tilde{Q}_\nb) - b_\ep (d)^{-2(r-1)}r_1^2
n_1^{-1}\sigma^2_{1,0,\epsilon} - b_\ep (d)^{-2(r-1)}r_2^2
n_2^{-1}\sigma^2_{0,1,\epsilon}  \notag \\
&=& K_{1,\mathbf{n}}b_\ep (d)^{-2(r-1)}n^{-1}\sigma^2_{1,0,\epsilon} + K_{2,%
\mathbf{n}}b_\ep (d)^{-2(r-1)}n^{-1}\sigma^2_{0,1,\epsilon}  \notag \\
&+& b_\ep (d)^{-2(r-1)}\sum_E \frac{\binom{r_1}{{l}}\binom{r_2}{{m}}\binom{%
n_1-r_1}{r_1-l}\binom{n_2-r_2}{r_2-{m}}} {\binom{n_1}{r_1}\binom{n_2}{r_2}}
\sigma^2_{{l},{m},\epsilon},
\end{eqnarray}
where $E:= \{(l,m): 0 \leq l \leq r_1, 0 \leq m \leq r_2, l+m \geq 2\}$, and
\begin{equation}
K_{1,\mathbf{n}} := r_1^2 p_\nb \left(\frac{\binom{n_1-r_1}{r_1-1}\binom{%
n_2-r_2}{r_2}}{\binom{n_1-1}{r_1-1}\binom{n_2}{r_2}} - 1 \right), \quad K_{2,%
\mathbf{n}} := r_2^2 (1-p_\nb) \left(\frac{\binom{n_1-r_1}{r_1}\binom{n_2-r_2%
}{r_2-1}}{\binom{n_1}{r_1}\binom{n_2-1}{r_2-1}} - 1 \right).  \notag
\end{equation}
Note that $K_{1,\mathbf{n}}, K_{2,\mathbf{n}}$ = $\mathrm{O}(n^{-1})$ $%
\mbox{ as } \mathbf{n} \to \infty$ so if $n\epsilon^d \to a$, $0<a \leq
\infty$, then \eqref{eq:bound}, \eqref{eq:varlim}, and \eqref{eq:varrem}
imply that $\mathrm{Var}(R_\nb) = \mathrm{O}((n^2\epsilon^d)^{-1})
\mbox{ as
} \mathbf{n} \to \infty$. In particular, for $a = \infty$,
\begin{equation}  \label{eq:ordo}
\mathrm{Var}(R_\nb) = \mbox{o}(n^{-1}) \Rightarrow n^{1/2} R_\nb \overset{%
\mathrm{P}}{\to} 0 \mbox{ as } \mathbf{n} \to \infty.
\end{equation}
By symmetry, we have from \eqref{cond}
\begin{eqnarray}  \label{eq:psi}
\psi_{1,0,\mathbf{n}}(x) = \frac{1}{r_1} \left(
p_{X,\epsilon}(x)^{r_1-1}p_{Y,\epsilon}(x)^{r_2} +(r_1-1)\mathrm{E}%
(\psi_\nb^{(2)}(x,X_2,\ldots,X_{r_1};Y_1,\ldots,Y_{r_2})) \right ).
\end{eqnarray}
Let $x$ be a continuity point of $p_X(x)$ and $p_Y(x)$. Then, changing
variables $y = x +\epsilon u$ and the bounded convergence theorem give
\begin{align}  \label{eq:psitwo}
\mathrm{E}(\psi_\nb^{(2)}(x,X_2&,\ldots,X_{r_1};Y_1,\ldots,Y_{r_2}) =
\mathrm{E}(\mathrm{E}(\psi_\nb^{(2)}(x,X_2,\ldots,X_{r_1};Y_1,%
\ldots,Y_{r_2})|X_2))  \notag \\
&= \int_{ R^d} I(d(x,y)\leq
\epsilon)p_{X,\epsilon}(y)^{r_1-2}p_{Y,\epsilon}(y)^{r_2}p_X(y)dy  \notag \\
&= \epsilon^d\int_{ R^d} I(d(0,u)\leq 1)p_{X,\epsilon}(x+\epsilon
u)^{r_1-2}p_{Y,\epsilon}(x+\epsilon u)^{r_2}p_X(x+\epsilon u)du \\
&\sim b_\ep (d)^{r-1} p_X(x)^{r_1-1} p_Y(x)^{r_2} \mbox{ as } \mathbf{n} \to
\infty.  \notag
\end{align}
From \eqref{eq:psi} we get that
\begin{equation}
\psi_{1,0,\mathbf{n}}(x) \sim b_\ep (d)^{r-1} p_X(x)^{r_1-1}p_Y(x)^{r_1} %
\mbox{ as } \mathbf{n} \to \infty,  \notag
\end{equation}
and hence
\begin{equation}  \label{eq:hlim1}
\lim_{\mathbf{n} \to \infty}h_\nb^{(1,0)}(x) = p_X(x)^{r_1-1}p_Y(x)^{r_2} -
q_\rb,
\end{equation}
and similarly,
\begin{equation}  \label{eq:hlim2}
\lim_{n \to \infty}h_\nb^{(0,1)}(x) = p_X(x)^{r_1}p_Y(x)^{r_2-1} - q_\rb.
\end{equation}
Let $\max(p_X(x),p_Y(x)) \leq C, x\in R^d$. Then $\max(p_{X,%
\epsilon}(x),p_{Y,\epsilon}(x)) \leq b_\ep (d) C, x \in R^d$. It follows
from \eqref{eq:psi} and \eqref{eq:psitwo} that $\psi_{1,0,\mathbf{n}}(x)
\leq b_\ep (d)^{r-1}C^{r-1}, x\in R^d$, and hence $h_\nb^{(1,0)}(x) \leq
2C^{r-1}$, $x \in R^d$. Similarly, we have that $h_\nb^{(0,1)}(x) \leq
2C^{r-1}$, $x \in R^d$. Therefore, $h_\nb^{(1,0)}(X_1)$ and $%
h_\nb^{(0,1)}(Y_1)$ are bounded random variables. Hence, from %
\eqref{eq:hlim1}, \eqref{eq:hlim2}, and the bounded convergence theorem we
obtain
\begin{eqnarray}
\mathrm{Var}{(h_\nb^{(1,0)}}(X_1)) \to \zeta_{1,0}, \quad \mathrm{Var}{%
(h_\nb^{(0,1)}}(Y_1)) \to \zeta_{0,1} \mbox{ as } \mathbf{n} \to \infty.
\notag
\end{eqnarray}
Let $Z_{\mathbf{n},i}:= n_1^{-1/2}h_\nb^{(1,0)}(X_i)$, $i =1, \ldots,n_1$,
and observe that, for $\delta > 0$,
\begin{equation}
\sum_{i=1}^{n_1} \mathrm{E} Z_{\mathbf{n},i}^2 = \mathrm{Var}%
(h_\nb^{(1,0)}(X_1)) \to \zeta_{1,0} > 0 \mbox{ as } \mathbf{n} \to \infty,
\notag
\end{equation}
\begin{eqnarray}
\lim_{\mathbf{n} \to \infty} \sum_{i=1}^{n_1} \mathrm{E}(|Z_{\mathbf{n}%
,i}|^2 I( |Z_{\mathbf{n},i}|>\delta)) & = & \lim_{\mathbf{n} \to \infty}
\mathrm{E}(|h_\nb^{(1,0)}(X_1)|^2I(|h_\nb^{(1,0)}(X_1)| >\delta n_1^{1/2}))
\notag \\
& \leq & \lim_{\mathbf{n} \to \infty} 4C^{2(r-1)}\mathrm{E}%
(I(|h_\nb^{(1,0)}(X_1)| >\delta n_1^{1/2})) = 0,  \notag
\end{eqnarray}
where the last equality follows from the boundedness of $h_\nb^{(1,0)}(X_1)$%
. The Lindeberg-Feller Theorem (see, e.g., Theorem 4.6, Durrett, 1991) gives
that
\begin{equation}
Z_{\mathbf{n},1}+\ldots+Z_{\mathbf{n},n_1} = n_1^{1/2} H_\nb^{(1,0)} \overset%
{\mathrm{D}}{\to} N(0,\zeta_{1,0}) \mbox{ as } \mathbf{n} \to \infty,  \notag
\end{equation}
and similarly $n_2^{1/2} H_\nb^{(0,1)} \overset{\mathrm{D}}{\to}
N(0,\zeta_{0,1}) \mbox{ as } \mathbf{n} \to \infty$. Hence, by independence
we get that
\begin{align}
&n^{1/2} (r_1 H_\nb^{(1,0)} + r_2 H_\nb^{(0,1)})  \notag \\
&= \frac{r_1}{p_\nb^{1/2}} n_1^{1/2} H_\nb^{(1,0)} + \frac{r_2}{%
(1-p_\nb)^{1/2}} n_2^{1/2} H_\nb^{(0,1)} \overset{\mathrm{D}}{\to} N(0,
\kappa) \mbox{ as } \mathbf{n} \to \infty,  \notag
\end{align}
so from \eqref{eq:ordo} and Slutsky's theorem,
\begin{align}
&n^{1/2}(\tilde{Q}_\nb-\tilde{q}_{\mathbf{r},\epsilon})  \notag \\
&= n^{1/2} (r_1 H_\nb^{(1,0)} + r_2 H_\nb^{(0,1)}) + n^{1/2} R_\nb \overset{%
\mathrm{D}}{\to} N(0,\kappa) \mbox{ as } \mathbf{n} \to \infty.  \notag
\end{align}
This completes the proof. \hfill $\Box$\newline
\medskip \noindent \textit{Proof of Theorem \ref{th:AsBias}: } The proof is
similar to that of the corresponding result in Leonenko and Seleznjev (2010)
so we give the main steps only. First we evaluate the bias term $B_\nb :=
\tilde{q}_{\mathbf{r},\epsilon} - q_\rb$. Let $V:=(V_1,\ldots, V_d)^{\prime}$
be an auxiliary random vector uniformly distributed in the unit ball $B_1(0)$%
, say, $V \in U(B_1(0))$. Then by definition, we have
\begin{eqnarray*}
B_\nb = \int_{R^d}p_{X,\epsilon}(x)^{r_1-1}p_{Y,\epsilon}(x)^{r_2}p_X(x)dx -
\int_{R^d} p_X(x)^{r_1} p_Y(x)^{r_2} dx = \mathrm{E}(D_\ep(X)),
\end{eqnarray*}
where
\begin{eqnarray*}
D_\ep(x)&:=& p_{X,\epsilon}(x)^{r_1-1}p_{Y,\epsilon}(x)^{r_2}-
p_X(x)^{r_1-1} p_Y(x)^{r_2} \\
& = & p_{X,\epsilon}(x)^{r_1-1}(p_{Y,\epsilon}(x)^{r_2}-p_Y(x)^{r_2}) +
p_Y(x)^{r_2}( p_{X,\epsilon}(x)^{r_1-1}-p_X(x)^{r_1-1}).
\end{eqnarray*}
It follows by definition that
\begin{eqnarray*}
D_\ep(x)& = &
P_1(x)(p_{Y,\epsilon}(x)-p_Y(x))+P_2(x)(p_{X,\epsilon}(x)-p_X(x)) \\
& = & \mathrm{E}\bigl(P_1(x)(p_{Y}(x-\epsilon
V)-p_Y(x))+P_2(x)(p_{X}(x-\epsilon V)-p_X(x))\bigr)
\end{eqnarray*}
where $P_1(x)$ and $P_2(x)$ are polynomials in $p_X(x),p_Y(x),\mathrm{E}%
(p_X(x - \epsilon V))$, and $\mathrm{E}(p_Y(x - \epsilon V))$. Now the
boundedness of $p_X(x)$ and $p_Y(x)$ and the H\"{o}lder condition for the
continuous differentiable cases imply
\begin{eqnarray*}
|D_\ep(x)| &\leq & C C_1 \epsilon^\al, C_1>0,  \notag
\end{eqnarray*}
and the assertion (i) follows.

For $\epsilon\sim c n^{-1/(2\alpha+d(1-1/r))}, 0<c<\infty$, $\alpha<d/2$, by
(i) and Theorem \ref{main}, we have
\begin{equation*}
B_\nb^2+v_\nb^2=\mathrm{O}(n^{-2\alpha/(2\alpha+d(1-1/r))}).
\end{equation*}
Now for some $C>0$ and any $A>0$ and large enough $n_1,n_2$, we obtain
\begin{equation*}
P\left(|\tilde{Q}_\nb-q_\rb|>A n^{-\alpha/(2\alpha+d(1-1/r))}\right)\le
n^{-2\alpha/(2\alpha+d(1-1/r))} \frac{B_\nb^2+v_\nb^2}{A^2}\le \frac{C}{A^2},
\end{equation*}
and the assertion (ii) follows. Similarly for $\alpha=d/2$.

Finally, for $\alpha > d/2$ and $\epsilon \sim L(n)n^{-1/d}$ and $%
n\epsilon^d \to \infty$, the assertion (iii) follows from Theorem \ref%
{thm:main} and the Slutsky theorem. This completes the proof. \hfill $\Box$

\subsection*{Acknowledgment}

The third author is partly supported by the Swedish Research Council grant
2009-4489 and the project "Digital Zoo" funded by the European Regional
Development Fund. The second author is partly supported of the Commissions
the European Communities grant PIRSES-GA 2008-230804 "Marie Curie Actions".


\bigskip {\noindent {\large \textbf{References}}} \medskip

\begin{reflist}

Baryshnikov, Yu., Penrose, M., Yukich, J.E.: Gaussian limits for generalized
spacings. Ann.\ Appl.\ Probab., {\bfseries 19} (2009) 158--185

Beirlant, J., Dudewicz, E.J., Gyorfi, L., van der Meulen, E.C.:
Non-parametric entropy estimation: An overview. Internat.\ Jour.\ Math.\ Statist.\
Sci. {\bfseries 6} (1997) 17--39

Copas, J.B., Hilton, F.J.: Record linkage: statistical models for matching
computer records. Jour.\ Royal Stat.\ Soc.\ Ser A. {\bfseries 153} (1990)
287--320

Costa, J., Hero, A., Vignat, C.: On Solutions to Multivariate Maximum $%
\alpha $-entropy Problems. Lecture Notes in Computer Science {\bfseries 2683}
(2003) 211--228

Demetrovics, J., Katona, G.O.H., Mikl\'os, D., Seleznjev, O., Thalheim, B.:
The average length of keys and functional dependencies in (random)
databases. In: {Proc.\ ICDT95} Eds: G.\ Gottlob and M.\ Vardi, LN in Comp.\
Sc.\ Springer: Berlin {\bfseries 893} (1995) 266--279

Demetrovics, J., Katona, G.O.H., Mikl\'os, D., Seleznjev, O., Thalheim, B.:
Asymptotic properties of keys and functional dependencies in random
databases. Theor.\ Computer Science {\bfseries 190} (1998) 151--166

Demetrovics, J., Katona, G.O.H., Mikl\'os, D., Seleznjev, O., Thalheim, B.:
Functional dependencies in random databases. Studia Scien.\ Math.\ Hungarica {%
\bfseries 34} (1998) 127--140

Durrett, R.: Probability: Theory and Examples. Brooks/Cole Publishing
Company, New York (1991)

Kapur, J.N.: Maximum-entropy Models in Science and Engineering. Wiley New
York (1989)

Kapur, J.N., Kesavan, H.K.: Entropy Optimization Principles with
Applications. Academic Press, New York (1992)


Kiefer, M., Bernstein, A., Lewis, Ph.\ M.\ Database Systems: An
Application-Oriented Approach. Addison Wesley (2005)

Koroljuk, V.S., Borovskich, Yu.V.: Theory of $U$-statistics. Kluwer, London
(1994)

Kozachenko, L.F., Leonenko, N.N.: On statistical estimation of entropy of
random vector. Problems Infor.\ Transmiss. {\bfseries 23} (1987) 95--101

Lee, A.J.: $U$-Statistics: Theory and Practice. Marcel Dekker, New York
(1990)

Leonenko, N., Pronzato, L., Savani, V.: A class of R\'{e}nyi information
estimators for multidimensional densities. Ann.\ Stat. {\bfseries36} (2008)
2153--2182. Corrections, Ann.\ Stat., 38, N6 (2010) 3837-3838

Leonenko, N., Seleznjev, O.: Statistical inference for the $\epsilon$%
-entropy and the quadratic R\'{e}nyi entropy. Jour.\  Multivariate
Analysis {\bfseries 101} (2010) 1981--1994


Neemuchwala, H., Hero, A., Carson, P.: Image matching using alpha-entropy
measures and entropic graphs. Signal Processing {\bfseries 85} (2005)
277--296

Pardo, L.: Statistical Inference Based on Divergence Measures. Chapman and
Hall (2006)

R\'enyi, A.: On measures of entropy and information. In: Proc. 4th Berkeley
Symp.\ Math.\ Statist.\ Prob.\ Volume 1. (1961)

R\'enyi, A.: Probability Theory. North-Holland, London (1970)

Seleznjev, O., Thalheim, B.: Average case analysis in database problems.
Methodol.\ Comput.\ Appl.\ Prob. {\bfseries 5} (2003) 395-418

Seleznjev, O., Thalheim, B.: Random databases with approximate record
matching. Methodol.\ Comput.\ Appl.\ Prob. {\bfseries 12} (2008) 63--89

Shannon, C.E.: A mathematical theory of communication. Bell Syst.\ Tech.\ Jour.\ {%
\bfseries 27} (1948) 379--423, 623--656

Szpankowski, W.: Average Case Analysis Of Algorithms On Sequences. John
Wiley, New York (2001)

Thalheim, B.: Entity-Relationship Modeling. Foundations of Database
Technology. Springer Verlag, Berlin, (2000)

Tsybakov, A.B., Van der Meulen, E.C.: Root-$n$ consistent estimators of
entropy for densities with unbounded support. Scandinavian Jour.\ Statistics {%
\bfseries 23} (1996) 75-83
\end{reflist}

\end{document}